\documentclass[12pt]{article}

\usepackage{amsmath}
\usepackage{amsthm}
\usepackage{amssymb}
\usepackage{amscd}
\usepackage[all]{xy}

\title{Calabi--Yau threefolds with infinitely many divisorial 
contractions} 
\author{Hokuto Uehara\footnote{Research Fellow of the Japan Society 
for the Promotion of Science.}}
\date{\empty}

\theoremstyle{plain}
\newtheorem{prop}{Proposition}[section]
\newtheorem{lem}[prop]{Lemma}

\newtheorem{thm}[prop]{Theorem}
\newtheorem{cor}[prop]{Corollary}
\newtheorem{con}[prop]{Conjecture}
\newtheorem{cla}[prop]{Claim}

\theoremstyle{definition}
\newtheorem{defn}[prop]{Definition} 
\newtheorem{prob}[prop]{Problem}

\theoremstyle{remark}
\newtheorem{rem}[prop]{Remark}
\newtheorem{ex}[prop]{Example}

\newcommand{\Int}{\operatorname{Int}}
\newcommand{\Sing}{\operatorname{Sing}}
\newcommand{\Aut}{\operatorname{Aut}}
\newcommand{\Bir}{\operatorname{Bir}}
\newcommand{\Pic}{\operatorname{Pic}}
\newcommand{\Exc}{\operatorname{Exc}}
\newcommand{\diag}{\operatorname{diag}}
\newcommand{\Q}{\mathbb Q}
\newcommand{\R}{\mathbb R}
\newcommand{\C}{\mathbb C}
\newcommand{\Z}{\mathbb Z}
\newcommand{\PP}{\mathbb P}
\newcommand{\sAbar}{\overline{\mathcal{A}}}
\newcommand{\NEbar}{\overline{NE}}
\newcommand{\Span}[1]{\left<#1\right>}

\begin{document}

\maketitle

\begin{abstract}
We study Calabi--Yau 3-folds with infinitely many divisorial contractions. We also 
suggest a method to describe Calabi--Yau 3-folds with the infinite automorphism group. 
\end{abstract}


\setcounter{section}{-1}
\section{Introduction}

\indent
 A smooth complex projective $n$-dimensional variety $X$ is a Calabi--Yau $n$-fold (C--Y $n$-fold) if $K_X=0$ and $h^1(\mathcal{O}_X)=0$.
If the Abundance Conjecture and the Minimal Model Conjecture are true, a $\Q$-factorial terminal $n$-fold $Y$ with Kodaira dimension
$\kappa (Y)=0$ is always birationally equivalent to a $\Q$-factorial terminal $n$-fold $X$ with $K_X\equiv 0$ (\cite{KMM}, \cite{KM}).
We can regard C--Y $n$-folds as special cases of this. As is well-known, for a smooth K3 surface $S$, the nef cone $\sAbar (S)$ is
rational polyhedral if and only if $\Aut S$ is finite (\cite{Sterk}). Moreover if a K3 surface $S$ with infinite $\Aut S$ contains 
a $-2$-curve, then $S$ contains infinitely many $-2$-curves (\cite{Kov}). In the same way, the Morrison Cone Conjecture (\ref{conj}) 
states that for a C--Y 3-fold $X$ the nef cone $\sAbar (X)$ is rational polyhedral if and only if $\Aut X$ is finite. 
By analogy with K3 surfaces and C--Y 3-folds, if a C--Y 3-fold $X$ with infinite $\Aut X$ admits a divisorial contraction,
it is highly likely that it admits infinitely many such. In addition to this, a C--Y 3-fold always admits a birational contraction
when its Picard number is more than 13 (\cite{HW}). In this context, it seems worthwhile to study C--Y 3-folds with infinitely many 
divisorial contractions. One of the aim of this article is to give a characterization of C--Y 3-folds which admit infinitely many divisorial
contractions (see Theorem \ref{structure 0}. See also Theorem \ref{structure} and Remark \ref{converse} for the precise statement). 
 
 Another aim of this article is to suggest a method to describe C--Y 3-folds $X$ with infinite $\Aut X$. 
If we have such $X$, then $\sAbar (X)\cap c_2 ^{\perp}\ne \{0\}$ (Remark \ref{rational polyhedral}), where $c_2(=c_2(X))$ is the second 
Chern class of $X$. If $\sAbar (X)\cap c_2 ^{\perp}$ contains the class of a rational divisor, it is likely (cf. Conjecture \ref{SC}) 
that some multiple of the divisor determines a nontrivial contraction $\varphi :X\to Y$ satisfying $\varphi ^*H\cdot c_2=0$
for an ample divisor $H$ on $Y$. We call such a contraction \emph{$c_2$-contraction}. 
In this context our first task to describe C--Y 3-folds with infinite $\Aut X$ is to:

\begin{list}{(i)}{}
\item describe C--Y 3-folds $X$ with infinite $\Aut X$ such that $X$ does not admit any nontrivial $c_2$-contractions. 
\end{list}

\noindent
I guess such $X$ has the small Picard number greater than 2. Secondly we should:

\begin{list}{(ii)}{}
\item classify C--Y 3-folds which admit a nontrivial $c_2$-contraction. 
\end{list}
 
\noindent 
Presumably we can do this because we have the remarkable classification of C--Y 3-folds $X$ admitting a $c_2$-contraction 
$\varphi :X\to Y$ in the case $\dim Y\ge 2$ by K. Oguiso (cf. \cite{OS} or Theorem \ref{O-S}). Next we should: 

\begin{list}{(iii)}{}
\item determine which C-Y 3-folds in the list obtained by (ii) have infinite $\Aut X$. 
\end{list}

\noindent
If we carry out these, we can describe all C--Y 3-folds with infinite $\Aut X$.
 
 In Section \ref{divisorial contractions}, we prove several lemmas for the latter use.
Let $\tilde I(= \tilde I _X)$ be the index of the set $\{\varphi _i\}_{i\in \tilde{I}}$ of all possible divisorial contractions on a 
C--Y 3-fold $X$ and let us denote the exceptional divisor of $\varphi _i$ by $E_i$. 
The most important lemma in Section \ref{divisorial contractions} is:

%
%

\begin{lem}\label{interesting 0}(= Proposition \ref{interesting} + Remark \ref{type II}.) 
Let $J$ be an infinite subset of $\tilde I$. Then there exist $1,2,3\in J$ such that $E_1+E_2+E_3$ is nef.
\end{lem}

\noindent
We use this lemma in Section \ref{section structure} to construct a nontrivial $c_2$-contraction on C--Y 3-folds with infinitely many 
divisorial contractions.

 In Section \ref{evidence}, we give a partial answer to the following conjecture.
Put $\sAbar (X)_{\epsilon} :=\bigl\{ x\in  \sAbar(X)\bigm| c_2 \cdot x\ge\epsilon H^2\cdot x \bigr\}$
for an ample divisor $H$ on $X$ and let $\epsilon$ be a positive real number.

%

\begin{con}\label{c_2 conjecture 0}(=Conjecture \ref{c_2 conjecture}.) 
Let $X$ be a C--Y 3-fold.
\renewcommand{\labelenumi}{(\roman{enumi})}
\begin{enumerate}
\item Let $\varphi :X\to Y$ be a contraction such that $\varphi ^*\sAbar (Y) \subset\sAbar (X)_{\epsilon}$. 
Then the cardinality of the set of such $\varphi$ is finite. 
\item Let $\varphi :X\to Y$ be a contraction such that $\varphi ^*\sAbar (Y) \subset\sAbar (X)_{\epsilon}$. 
Then $\sAbar (Y)$ is rational polyhedral.
\end{enumerate}
\end{con}

\noindent
If $\Aut X$ is infinite, then $\sAbar (X)$ is not rational polyhedral (Remark \ref{rational polyhedral}).
Hence Conjecture \ref{c_2 conjecture 0} means the shape of $\sAbar (X)$ is complicated near $\sAbar (X)\cap c_2^{\perp}$. 
We expect this ``complexity'' produces a rational point on $\sAbar (X)\cap c_2^{\perp}\backslash  \{0\}$.     
 
 In Section \ref{section structure}, we consider C--Y 3-folds with infinitely many divisorial contractions.
Define $\tilde I_{c_2 *0}:=\bigl\{i\in \tilde I\bigm| E_i \cdot c_2 *0 \bigr\}$, where $*$ is $<,=$ or $>$. 
The main result of Section \ref{section structure} is:

 %
 %

\begin{thm}\label{structure 0}(See Theorem \ref{structure} for the precise statement.)
Assume that $\tilde{I}_{c_2=0}$ is infinite for a C--Y 3-fold $X$. Then there exist a K3 surface $S$
containing infinitely many smooth rational curves, an elliptic curve $E$ and a finite Gorenstein automorphism group $G$ of $S\times E$ such that
$X$ is birational to $(S\times E)/G$. 
\end{thm}

\noindent
 In the proof of Theorem \ref{structure 0} we use Lemma \ref{interesting 0} to prove the existence of a nontrivial $c_2$-contraction on $X$ 
and we use the Oguiso's classification to determine the structure of $X$. Hence Theorem \ref{structure 0} is regarded as a realization of 
the method to describe C--Y 3-folds with infinite $\Aut X$ we mention above. 
 
 Finally, in Section \ref{construction} we construct C-Y 3-folds with $|\tilde{I}_{c_2=0}|=\infty$. In passing, 
we show that the set $\tilde{I}_{c_2<0}$ is always finite in Corollary \ref{fano} and Remark \ref{type II}. I do not know any examples
of C-Y 3-folds with $|\tilde{I}_{c_2>0}|=\infty$.

\paragraph{Notation and Convention}

\noindent
\begin{enumerate}
\renewcommand{\labelenumi}{(\roman{enumi})}

\item  When a normal projective variety $X$ over $\C$ has
at most rational Gorenstein singularities and it satisfies $h^1(\mathcal{O}_X)=0$ and $K_X=0$, we call it a C--Y model. $X$ always means a C--Y 3-fold and a C--Y model means a 3-dimensional C--Y model throughout this paper unless we specify
otherwise.

\item For a $n$-dimensional projective variety $X$, let $\mathcal{A} (X)$
denote the cone generated by ample divisors in $N^1 (X)$ and
$\mathcal{A}^e(X)$ denotes the effective nef cone, namely, the cone
generated by nef effective divisors in $N^1(X)$. Let us denote the cone $\bigl\{ x\in N^1 (X)\bigm| x^n =0\bigr\}$ by $\mathcal{W}$. Suppose the symbol $*$ denotes $>,\ge$ etc. For a real divisor $D$ on $X$ and a constant $c$, set $D_{* c}:=\bigl\{ z\in N_1(X)\bigm| (D\cdot z) * c\bigr\}\cup\{0\}$. Moreover $[D]$ denotes the element in $N^1(X)$ corresponding to $D$. For a real 1-cycle $z$, define the subspace $z_{* c}$ of $N^1(X)$ and the class $[z]\in N_1(X)$ in the similar way. Define $\overline{NE}(X)_{D*0}:= \overline{NE} (X)\cap D_{*0}$.    

\item For a C--Y 3-fold $X$, we can regard the second Chern class $c_2(X)$
as a linear form on $H^2(X,\Z)$. We often abbreviate it by $c_2$ in this
article. As is well-known, $c_2\cdot x\ge 0$ for all $x\in \sAbar(X)$ by 
Y. Miyaoka (\cite{Miyaoka}). We define $\sAbar (X)_{\epsilon} :
=\sAbar(X)\cap (c_2 -\epsilon H^2)_{\ge 0}$ for a fixed ample
divisor $H$ and a positive real number $\epsilon$. 

\item We use the terminology \emph{terminal}, \emph{canonical}, \emph{klt}
(Kawamata log terminal), \emph{lc} (log canonical) and \emph{plt} (purely
log terminal) for a log pair $(X,\Delta)$ in the sense in \cite{KM}, but
we always assume that $\Delta$ is effective in these definitions. Klt is
same as log terminal in \cite{KMM}. We also use the terminology
\emph{semismooth} in the sense in \cite{FA}.   

\item The term \emph{contraction} means a surjective morphism between normal projective varieties with connected fibers and 
thus contractions consist of the fiber space case and the birational contraction case. Let $I_X(=I)$ be the index of the set $\{\varphi _i\colon X\to Y_i\}_{i\in I}$ of all possible birational contractions of type III on a C--Y 3-fold $X$ (see Definition \ref{I,II,III} for this terminology).
For $i\in I$, let $E_i$ be the exceptional divisor of $\varphi _i$, $C_i$
the irreducible curve $\varphi _i(E_i)$ and $F_i$ a general fiber of
$\varphi _i |_{E_i} \colon E_i \to C_i$. It is known that $E_i\cdot F_i =-2$. Furthermore let us denote by $V_i$ the image of the closed cone of curves $\overline {NE} (E_i)$ under the natural map $N_1(E_i)\to N_1(X)$. We know that $V_i$ is a 2-dimensional cone (see Fact (iii)) generated by the rays $\R _{\ge 0}[F_i]$ and $\R _{\ge 0}[v_i]$, where $v_ i$ is a real 1-cycle.

\item We denote the biregular (respectively, birational) automorphism group of a variety $X$ by $\Aut X$ (respectively, $\Bir X$).

\item If $V$ is given as $V_{\Q} \otimes \R$ for some $\Q$-vector space
$V_{\Q}$, a \emph{rational polyhedral} cone is a closed cone generated by
a finite set of rational points. A cone $\mathcal{C}$ is \emph{locally
rational polyhedral at a point} $x$ if there is a neighborhood $U$ of $x$
and a rational polyhedral cone $\mathcal{D}$ such that $\mathcal{C}\cap
U=\mathcal{D}\cap U$. Let $\mathcal{E}$ be a open cone in $V$. We say that
a cone $\mathcal{C}$ is \emph{locally rational polyhedral in
$\mathcal{E}$} if $\mathcal{C}$ is a rational polyhedral at every point in
$\mathcal{E}$. 
\end{enumerate}

\noindent
\textit{Acknowledgment.} 
Part of this paper was written while I was at the University of Warwick from April 2000 to March 2001 and had a short visit at the University of Bayreuth. I would like to thank both universities (in particular Professors M. Reid and Th. Peternell) for providing a pleasant environment.
I wish to thank Professors Y. Kawamata and Masa-hiko Saito for pertinent comments. I am especially grateful to Professors K. Oguiso and B. Szendr\"oi for their generosity with their time and expertise. Last but not least, I would like to express gratitude to my friends, K. Inui, M. Fukushima and D. Ryder for their warm encouragement.
\\


\section{Divisorial contractions on C--Y 3-folds}\label{divisorial contractions}
 We say that a birational contraction $\varphi \colon X\to Y$ between normal projective varieties is \emph{primitive} if $\rho (X/Y)=1$. We classify a primitive birational contraction on a $\mathbb{Q}$-factorial C--Y model according to the dimensions of its exceptional set and its image. 


\begin{defn}\label{I,II,III}
We say that a primitive birational contraction on a (3-dimensional) C--Y model is \emph{of type I} if it contracts only finitely many curves, \emph{of type II} if it contracts an irreducible surface to a single point and \emph{of type III} if it contracts an irreducible surface to a curve.
Hence a primitive birational contraction is, so called, a small (respectively, divisorial) contraction if it is of type I (respectively, type II or III). Every birational contractions on a $\mathbb{Q}$-factorial C--Y model is one of types I, II and III.
\end{defn}

 Let $\varphi \colon X\to Y$ be a birational contraction on a $n$-dimensional C--Y model $X$. Let $H$, $H'$ denote ample divisors on $X$,
$Y$ respectively. Since $\Delta:=-H+m\varphi ^* H'$ is effective for sufficiently large $m$, the pair $(X, \epsilon \Delta)$ defines a log 
variety with klt singularities for $0<\epsilon \ll 1$. Therefore we can regard $\varphi$ as a $K_X +\epsilon \Delta$-extremal face contraction
and so we may apply theory of the log Minimal Model Program (log MMP) to study $\varphi$. All of the following facts come from theory of the
log MMP (\cite{KMM}, \cite{KM}).  

\paragraph{Fact}

\noindent
\begin{enumerate}
\renewcommand{\labelenumi}{(\roman{enumi})}
\item
Since $-(K_X +\epsilon \Delta)$ is $\varphi$-ample, the cone $\overline{NE}(X/Y)$ is rational polyhedral by the cone theorem.

\item
Since every extremal \emph{face} contraction can be decomposed into extremal \emph{ray} contractions, we can write $\varphi =\psi_m \circ
\cdots \circ \psi _1$, where $\psi _i$ is a primitive contraction and
$m=\rho (X/Y)$. A contraction $\varphi$ corresponds to a codimension $m$
face
$\Delta _m$ of $\sAbar (X)$, not entirely contained in $\mathcal{W}$,
which is just the image of $\sAbar (Y)$ under the injection $\varphi ^*
\colon N^1 (Y)\to N^1 (X)$. Thus a decomposition of $\varphi$ corresponds
to a sequence of faces $\Delta _0 :=\sAbar (X) >\Delta_1>\cdots >\Delta _m$,
where $\Delta _i$ is a codimension 1 face of $\Delta _{i+1}$.

\item
Since the image of $\varphi ^* \colon\Pic (Y)\to \Pic(X)$ coincides with 
\[\bigl\{ D\in \Pic (X)\bigm| D \cdot z=0 \ \textup{for all} \ z \in (\varphi ^* {H'})^{\perp}\cap \overline{NE}(X)\bigr\}\] 
and since $X$ is a C--Y model, $Y$ is also a C--Y model. We also obtain an exact sequence 
\[0\to N_1(X/Y)\to N_1(X)\to N_1(Y)\to 0.\]         
Assume that $\dim X=3$. Pick $i\in I$. From the exact sequence above, we know that $V_i$ is a 2-dimensional cone in $N_1 (X)$. 

\item
Let $X$ be a C--Y 3-fold and $L$ an \emph{effective} nef divisor on it. Since $(X, \epsilon L)$ is a klt pair for $0<\epsilon \ll 1$ and $K_X +\epsilon L$ is nef, we know that $L$ is semiample by the log abundance theorem (\cite{KeMM}, see also \cite{Oguiso93}). 

%
%
%
 
\begin{con}\label{SC} 
Let $X$ be a C--Y 3-fold and $L$ a nef divisor on it. Then $L$ is semiample.
\end{con}

\noindent If $L\cdot c_2 >0$, we can show that $L$ is effective (\cite{W94}). So in this case, Conjecture \ref{SC} is true. 

\item
By the cone theorem for klt pairs, the nef cone $\sAbar (X)$ is locally rational polyhedral inside the cone $\mathcal{W}$. See \cite{K88}, \cite{K97} and \cite{W94} for the proof. 
\end{enumerate}
  
 In passing, for a C--Y 3-fold $X$ and an effective divisor $\Delta$ on it
such that the pair $(X,\Delta)$ has at most klt singularities, if every
$K_X +\Delta$-extremal ray corresponds to a divisorial contraction, the
number of $K_X +\Delta$-extremal rays is finite by the observation in
Fact (iii). On the other hand, the pair of the C--Y 3-fold $X$ constructed
by C. Schoen (cf. \cite{Namikawa}) and some effective divisor $\Delta$ on
$X$ gives an example where $\overline{NE} (X)_{K_X +\Delta <0}$ contains
infinitely many extremal rays corresponding to contractions of type I
(\cite{Namikawa}). This supplies a negative answer for the problem stated in (4-2-5) \cite{KMM},
i.e. for a klt pair $(X,\Delta)$ with $\kappa (X,K_X +\Delta)\ge 0$, is the number 
of $K_X +\Delta$-extremal rays finite?
But I still feel (4-2-5) ibid. is affirmative when $\Delta$ is trivial.  
       
 We have the following result by V. V. Nikulin \cite{Nikulin}, p282.

\begin{prop} \label{special}
The sets $I^1:=\bigl\{ i\in I\bigm| E_i$ is an exceptional divisor of two different divisorial contractions $\bigr\}$ and $I^2:=\bigl\{ i\in I\bigm|$ there exists $j\in I$ such that \emph{either} $E_i\cdot F_j >0$ and $E_j\cdot F_i=0$ \emph{or} $E_j\cdot F_i >0$ and $E_i\cdot F_j=0\bigr\} $ are finite. 
\end{prop}
 
%
%
%

\begin{lem}\label{fundamental} 
Let $X$ be a $\mathbb{Q}$-factorial C--Y model with its Picard number $\rho$.
Define $K_i:=\bigl\{ j\in I\bigm| E_i \cap E_j
\neq \emptyset\bigr\}$ for $i\in I$.
\renewcommand{\labelenumi}{(\roman{enumi})}
\begin{enumerate}
\item 
Assume $J\subset I$. If $|J|\geq \rho$, there exist
$i$, $j\in J$ such that $E_i\cap E_j$ is not empty. 
\item
There is no subset $J \subset I$ such that $J$ satisfies the
following property $(\ast)$.

 $(\ast)$ Assume that we have $1,\dots ,n \in J$ such that $i\in J\backslash
\bigcup_{k=1}^{i-1} K_k$ for all $i\leq n$. Then $J\backslash \bigcup_{k=1}^{n} 
K_k \neq \emptyset$.
\item
 Assume $J \subset I$ such that $|J|=\infty$. Then there exists $i\in J$ such that $|K_i \cap J|=\infty$. In particular, there exists an infinite subset $J'\subset J$ such that $E_i \cap E_j$ is not empty for all $i,j \in J'$.
\end{enumerate}
\end{lem}

\begin{proof}
(i) Assume that we have elements $1,\dots ,\rho\in J$ such that  $E_i \cap E_j$ is empty for all $i\ne j$.  Then there exists a nontrivial relation 
$\Sigma_{k=1}^{\rho}a_k E_k+ a_0 H\equiv 0$ 
for $a_k\in\R$ and some ample divisor $H$. Then because $E_i\cdot F_j =0$ if and only if $i\ne j$, the numbers $a_k\cdot a_0>0$ for all $k$. This is absurd, since $(\Sigma a_k E_k+ a_0 H)\cdot H^2\ne 0$.

(ii) If $J$ satisfies ($\ast$) then we have $1,\dots,\rho \in J$ such that
$k\notin \bigcup_{i=1}^{k-1} K_i$ for all $k \leq \rho$. This contradicts (i).

(iii) Assume that $K_i\cap J$ is finite for all $i\in J$. By $|J|=\infty$, $J$ satisfies ($\ast$) in (ii). The second statement follows from the first one.
\end{proof}

\begin{rem}\label{type II}
Every exceptional divisor of a birational contraction of type II does not meet each other. Therefore the number of contractions of type II is finite by the same proof of (i) above.
\end{rem}
 
\begin{lem}\label{special2}
For general $i\in I$ (namely, all but a finite number of $i\in I$) $\overline{NE}(X)=\overline{NE}(X)_{E_i\ge 0}+\R _{\ge 0}[F_i] $. 
\end{lem}

\begin{proof}
It is enough to check the finiteness of $J:=I \backslash (I^1 \cup I^2 \cup \bigl\{ i\in I\bigm| \overline{NE}(X)=\overline{NE}(X)_{E_i\ge 0}+\R _{\ge 0}[F_i] \bigr\} )$. For $i\in J$, not only $\R _{\ge 0}[F_i]$ but also $\R _{\ge 0}[v_i]$ is a $K_X +\epsilon E_i$-extremal ray. Then $\R _{\ge 0}[v_i]$ determines a birational contraction of type I. If $J$ is infinite, there exists an infinite subset $J'\subset J$ such that $E_i \cap E_j$ is not empty for all $i,j \in J'$ by Lemma \ref{fundamental}.
Then $\R _{\ge 0}[v_i] =\R _{\ge 0}[v_j]$ for all $i, j \in J'$. Let $\varphi \colon X \to Y$ be the associated contraction of type I and $H$ a general hyperplane section on $Y$, and define $l_i :=\varphi (E_i)|_H$ for $i\in J'$. Then since $l_i\cdot l_j =0$ on $H$ if and only if $i\ne j$, the $l_i$'s are linearly independent in $N_1(H)$. This is absurd.
\end{proof} 
Pick $i\in I$. Define $t_i=$min$\bigl\{ t\in \R \bigm| E_i +tH$ is nef $\bigr\}$, where $H$ is a fixed ample divisor on $X$. $\{t_i\}$ denotes the round up of $t_i$. 

%
%
%

\begin{lem}\label{nefvalue}
$t_i \leq 4$ for all $i \in I$.
\end{lem}

\begin{proof} 
If $E_i$ is normal, $E_i$ has at most RDP. By the inversion of adjunction, $(X,E_i)$ has at most plt singularities. If $E_i$ is non-normal, $E_i$ is semismooth
(\cite{W97?}). Then we use the inversion of adjunction again and know $(X,E_i)$ has at most lc singularities. In both cases,
we can apply the rationality theorem (\cite{KMM}) for the klt pairs $(X,
(1-\epsilon )E_i)$ for sufficiently small positive rational numbers
$\epsilon$ and we obtain the statement.
\end{proof}

%
%
%

\begin{lem}\label{compact}
Let $J\subset I$ and let $H$ be an ample divisor on $X$. Assume that there exist an integer $N$ and $z \in \overline {NE}(X)$ such that $z\cdot E_i \le N$ for all $i \in J$.
\renewcommand{\labelenumi}{(\roman{enumi})}
\begin{enumerate} 
\item Let $\epsilon$ be a positive real number. Then the set $J_\epsilon (z):=\bigl\{ i\in J\bigm| \varphi _i^* \sAbar (Y_i)\subset (z-\epsilon H^2)_{\ge0} \bigr\}$ is finite.

\item If $z$ is in the interior of $\overline {NE}(X)$, $J$ is finite. 
\end{enumerate}
\end{lem}

\begin{proof}
(i) By Lemma \ref{special2}, we may assume that $E_i+t_i H \in\varphi _i^* \sAbar (Y_i)$ for all $i\in J_\epsilon (z)$.
Then we get 
$(E_i +\{t_i\}H)\cdot (z-\epsilon H^2) \ge (\{t_i \}-t_i)H \cdot (z-\epsilon H^2) \ge 0$ 
and $(E_i +\{t_i\}H)\cdot z\le N+4H\cdot z =:c$. Thus $E_i +\{t_i\} H \in (z-\epsilon H^2)_{\ge 0} \cap z_{\le c} \cap \sAbar (X)$. Since $(z-\epsilon H^2)_{\ge 0} \cap z_{\le c} \cap \sAbar (X)$ is a compact set, $J_{\epsilon}$ is finite.

(ii) This is the special case of (i).
\end{proof}  

Let $D$ be a prime divisor on $X$. By the Serre duality for a Cohen-Macaulay surface $D$,
\[ \chi (\mathcal{O}_{D})=\chi (\omega _{D})
=\chi (\mathcal{O}_{D}(D))=\chi (\mathcal{O}_X(D)).\]
Combining this equality with the Riemann-Roch theorem for a C--Y 3-fold $X$, we obtain:
  
%
%
%

\begin{lem}\label{formula} 
For a prime divisor $D$ on $X$, we have
\[ \chi (\mathcal{O}_{D})=(1/6)D ^3+(1/12)D\cdot c_2.\] 
\end{lem}

\noindent The following proposition is a key to prove Theorem \ref{structure}.

\begin{prop}\label{interesting}
Let $J$ be an infinite subset of $I$. Then there exist $1,2,3\in J$ such that $E_1+E_2+E_3$ is nef.
\end{prop}

\begin{proof}
We may assume that $\overline{NE}(X)=\overline{NE}(X)_{E_i\ge 0}+\R _{\ge 0}[F_i] $ for all $i\in J$ by Lemma \ref{special2} and that $E_i \cdot F_j >0$ for all different $i,j\in J$ by Proposition \ref{special} and Lemma \ref{fundamental}(iii). Pick $1,2,3\in J$. Then $(E_1+E_2+E_3) \cdot F_i \ge0$ for $i=1,2,3$. Thus $E_1+E_2+E_3$ is nef.
\end{proof}

\noindent 
Note that the nef divisor $E_1+E_2+E_3$ is semiample by Fact (iv). By Proposition \ref{interesting}, the set $\bigl\{i\in I\bigm| E_i \cdot z <0 \bigr\}$ is finite for a pseudo-effective element $z \in N_1(X)$, i.e. $z\cdot x\ge 0$ for all $x\in \sAbar(X)$.  

%
%

\begin{cor}\label{fano}
The sets $I_{c_2<0}:=\bigl\{i\in I\bigm| E_i \cdot c_2 <0 \bigr\}$, $\bigl\{i\in I\bigm| E_i$ is a Hirzeburch surface $\bigr\}$ and $I_{dP}:=\bigl\{i\in I\bigm| E_i$ is a generalized del Pezzo surface $\bigr\}$ are finite.
\end{cor}

\begin{proof}
Because $c_2$ is pseudo-effective on minimal model 3-folds by \cite{Miyaoka}, the set $I_{c_2<0}$ is finite. For $i\in I$ such that $E_i$ is a Hirzeburch surface, $E_i \cdot c_2 =-4$ by Lemma \ref{formula}. Next suppose that $I_{dP}$ is infinite. By Proposition \ref{special} and Lemma \ref{fundamental}(iii), we may assume that $E_i \cdot F_j >0$ for all different $i,j\in I_{dP}$. Then there exists a real 1-cycle $v$ such that $\R_{\ge 0}[v]=\R_{\ge 0}[v_i]$ for all $i\in I_{dP}$. This is absurd, since $E_i \cdot v<0$ for all $i\in I_{dP}$.
\end{proof}


\section{The second Chern class and the nef cone}\label{evidence}
Let us remember the following conjecture of D. Morrison concerning the finiteness properties of the nef cones (\cite{Morrison}, \cite{K97}). We refer to \ref{conj} as the Morrison Cone Conjecture.


\begin{con}\label{conj}
Let $X$ be a C--Y n-fold. The number of the $\Aut X$-equivalence classes of faces of the effective nef cone $\mathcal{A}^e (X)$ corresponding to birational contractions or fiber space structures is finite. Moreover, there exists a rational polyhedral cone $\Pi$ which is a fundamental domain for the action of $\Aut X$ on $\mathcal {A}^e (X)$ in the sense that
\renewcommand{\labelenumi}{(\roman{enumi})}
\begin{enumerate}
\item $\mathcal{A} ^e (X)=\bigcup_{\alpha \in \Aut X} \alpha _* \Pi$,
\item $\Int \Pi \cap \alpha _* \Int \Pi =\emptyset$ unless $\alpha _*=id$.
\end{enumerate}
\end{con}

Let $H$ be a nef and big divisor on a (3-dimensional) C-Y model $Y$. 
Set $\Aut (Y,H):=\bigl\{\alpha\in \Aut Y\bigm|\alpha_* H\equiv H\bigr\}$. 

%
%
%

\begin{lem}\label{nef&big}
Let $Y$, $H$ be as above. Then the group $\Aut (Y,H)$ is finite.
\end{lem}

\begin{proof}
Let $\varphi\colon Y\to Z$ be the birational contraction defined
 by the free complete linear system $mH$ for sufficiently large integer $m$. Take an element of $\Aut (Y,H)$. Then it descends to an element of $\Aut (Z, H')$, where $H'$ is an ample divisor on $Z$ such that $\varphi ^*H'=mH$. On the other hand, the natural map $\Bir Y\to \Bir Z$ is injective, hence it is enough to prove the finiteness of $\Aut(Z,H')$. Grothendieck proved that $\Aut (Z,H')$ is a projective scheme, in particular, it has finitely many components. On the other hand, because $H^0(Y, T_Z)=0$ by Corollary 8.6, \cite{K85}, $\Aut Z$ is discrete and thus $\Aut (Z,H')$ is finite.  
\end{proof}

\begin{rem}\label{rational polyhedral}
If $c_2$ is positive on $\sAbar (X)\backslash \{0\}$ or if $\sAbar (X)$ is rational polyhedral, then since we can find an ample divisor $H$ such that $\Aut X=\Aut (X,H)$, $\Aut X$ is finite (\cite{W97}). Consequently if the Morrison Cone Conjecture is true for C--Y 3-folds $X$, $\sAbar (X)$ is rational polyhedral if and only if $\Aut X$ is finite.
\end{rem}

\noindent We study birational contractions of type III whose exceptional divisors are non-normal. If the Morrison Cone Conjecture is true,
we can bound the numbers $E_i^3$ and $E_i\cdot c_2$ for $i\in I$. In fact, for non-normal exceptional divisors $E_i$ we can prove
(without assuming the Morrison Cone Conjecture):

  %
  %
 
\begin{prop}\label{bounds}
$7-7h^{1,2}(X)\le E_i ^3\le 7$ and $-2\le E_i \cdot c_2 \le 6h^{1,2}(X)-2$ for all $i\in I$ such that $E_i$ is non-normal.
\end{prop}

\begin{proof} 
Fix $i\in I$ such that $E_i$ is non-normal and let $E$, $C$ denote $E_i$, $C_i$ respectively. Since $E$ is non-normal, $E$ is semismooth
and $C_0:=\Sing (E)$ is an irreducible smooth curve, which gives a section of $E\to C$ (\cite{W97?}). Let $\psi:Z\to X$, $E'$ and $D$ be the 
blowup along $C_0$, the strict transform of $E$ on $Z$ and the exceptional divisor of $\psi$ respectively. Let us also define $p:=\psi |_{E'}$ and
$C'_0:=p^{-1}(C_0)$ with the reduced structure. By local calculation, we can check easily that $p$ gives the normalization of $E$ and that $D$ and
$E'$ meet transversally, in particular, $D|_{E'}=C_0'$. Let $E'\to C'\to C$ be the Stein factorization of the morphism $E' \to E \to C$, then we
know that $E'$ is a $\mathbb{P}^1$-bundle over a smooth curve $C'$ and $C'\to C$ is a double cover. From these facts, we know that $C'_0$ is a
section of the $\mathbb{P}^1$-bundle. 
 
 Let $F$ be a ruling of the Hirzebruch surface $D$ over $C_0$. Because $\psi^* E|_D\cdot F=0$, $\psi^* E|_D$ is numerically proportional to 
$F$ on $D$ and so $0=(\psi^* E)^2\cdot D$. Then we have
\[0={E'}^2\cdot D+4E'\cdot D^2+4D^3.\]  
Furthermore because of $K_Z=D$ and the adjunction formula, we obtain
\[8(1-g(C'))=K_{E'}^2=D^2\cdot E'+2D\cdot {E'}^2+{E'}^3,\] 
\[2g(C')-2=(K_{E'}+C' _0)\cdot C' _0=2D^2\cdot {E'}+{E'}^2\cdot D\]
and 
\[8(1-g(C))=K_D^2=4D^3.\]
By these equalities, we get
\[E^3=(E'+2D)^3=7-3g(C')-4g(C).\]
By the fact that $g(C')\le h^{1,2}(X)$ (\cite{G97}), we get the bound of $E^3$.
On the other hand, because every fiber of $\varphi|_E\colon E\to C$ is a conic we have $R^i  \varphi_* \mathcal{O}_E =0$ for $i>0$. Thus we know $\chi (\mathcal{O}_E)=\chi (\mathcal{O}_C)$ and therefore
\[E\cdot c_2=12\chi (\mathcal{O}_E)-2E ^3=6g(C')-4g(C)-2\]
by Lemma \ref{formula}.
We use $g(C')\le h^{1,2}(X)$ again to obtain the bound of $E\cdot c_2$.
\end{proof}

\begin{rem}
We use the notation in the proof above. It seems worthwhile to restate the following formulae, that is, $E^3=7-3g(C')-4g(C)$ and
$E\cdot c_2=6g(C')-4g(C)-2$.
\end{rem}

%

\begin{con}\label{c_2 conjecture} (cf. Problem 3, \cite{W97})
\renewcommand{\labelenumi}{(\roman{enumi})}
\begin{enumerate}
\item Let $\varphi :X\to Y$ be a contraction such that $\varphi ^*\sAbar (Y) \subset\sAbar (X)_{\epsilon}$. 
Then the cardinality of the set of such $\varphi$ is finite. 
\item Let $\varphi :X\to Y$ be a contraction such that $\varphi ^*\sAbar (Y) \subset\sAbar (X)_{\epsilon}$. 
Then $\sAbar (Y)$ is rational polyhedral.
\end{enumerate}
\end{con}

\noindent
If $\Aut X$ is finite, the Morrison Cone Conjecture implies that the nef cone $\sAbar (X)$ is rational polyhedral. 
Hence obviously Conjecture \ref{c_2 conjecture} is true for such $X$ (modulo the Morrison Cone Conjecture).
If $\Aut X$ is infinite, then by Conjecture \ref{c_2 conjecture} we can expect the shape of the nef cone $\sAbar (X)$ 
is complicated near $\sAbar (X)\cap c_2 ^{\perp}$ (see also the argument after Problem \ref{existence of c_2}). 

 If we have a bound of the number $E_i\cdot c_2$ for $i\in I$, Conjecture \ref{c_2 conjecture}(i) is affirmative in the case when $\varphi$ is a birational contraction of type III, due to Lemma \ref{compact}(i).

%
%

\begin{thm}\label{ans. of c_2 conjecture}
Conjecture \ref{c_2 conjecture}(i) is affirmative in the following cases:
\renewcommand{\labelenumi}{(\roman{enumi})}
\begin{enumerate}
\item 
$\varphi$ is a fiber space (\cite{OP}).
\item 
$\varphi$ is a birational contraction of type III whose exceptional divisor is non-normal. 
\end{enumerate}
\end{thm}

  %
  %
  %
  %

\begin{thm}\label{ans. of c_2 comjecture II}
Conjecture \ref{c_2 conjecture}(ii) is affirmative in the following cases:
\renewcommand{\labelenumi}{(\roman{enumi})}
\begin{enumerate}
\item
$\varphi$ is a fiber space.
\item
Assume that the Morrison Cone Conjecture holds true and $\varphi$ is a birational contraction. 
\end{enumerate}
\end{thm}

\begin{proof}
(i) We may assume $\rho (Y)\ge 2$ so in particular $\dim Y=2$. By our assumption
and Theorem \ref{ans. of c_2 conjecture}(i) we know that $Y$ admits at most finitely many contractions. 
By Theorem 3.1 in \cite{Oguiso93} there exists a nonzero effective divisor $\Delta =\sum a_i D_i$ ($a_i >0$, $D_i$ a prime divisor) such 
that $(Y,\Delta)$ is a klt pair and $K_Y +\Delta\equiv 0$. 
Let $R=\R _{\ge 0}[z]$ be a \emph{geometrically extremal ray} of the cone $\NEbar (Y)$, where $z$ is a real 1-cycle 
(by the definition of a \emph{geometrically extremal ray}, if $z_1+z_2\in R$ for $z_1,z_2\in \NEbar (Y)$ we have $z_1,z_2\in R$.
Of course an extremal ray in the Minimal Model theory is geometrically extremal).
Note that $R$ is a $K_Y$-extremal ray if $K_Y\cdot z<0$, and $R$ is a $K_Y +\Delta +\epsilon D_i$-extremal ray for
some $i$ and $0< \epsilon \ll 1$ if $K_Y\cdot z >0$.
Now we prove that $\sAbar (Y)$ is rational polyhedral by the induction for $\rho (Y)$.
Denote the set of the geometrically extremal rays $R$ with $R\subset K_Y ^{\perp}$ by $\mathcal{S}$.
If $\mathcal{S}=\emptyset$ we have a contraction $f :Y\to Z$ for any geometrically extremal rays $R$ 
such that $f$ contracts only $R$.
So the proof is done by Theorem \ref{ans. of c_2 conjecture}(i). Hence we may assume $\mathcal{S}\ne\emptyset$. 
Pick $R(=\R _{\ge 0}[z])\in\mathcal{S}$. It is enough to show that we can take the real 1-cycle $z$ as a rational one and 
$\mathcal{S}$ is a finite set.
Since the cone $\NEbar (Y)$ is generated by the finitely many $K_Y$-extremal rays and the subcone $\NEbar (Y)_{K_Y \le 0}$, 
there exists a contraction $f(=f_R) :Y\to Z$ associated to a $K_Y$-extremal ray such that
$\R_{\ge 0}[z]+\R_{\ge 0}[F]=(f ^*L)^{\perp}\cap \NEbar (Y)$, where $F$ is a curve contracted by $f$ and $L$ is a nef $\R$-divisor on $Z$.
We can check that $f_*R$ is a geometrically extremal ray of the cone 
$\NEbar (Z)$ by using the exact sequence $0\to \Span{[F]} _{\R} \to N_1(Y)\to N_1(Z)\to 0$. Hence by the induction hypothesis 
(the finiteness of geometrically extremal rays of $\NEbar (Z)$),
there exists only finitely many $R_1\in\mathcal{S}$ such that $f_R=f_{R_1}$ (here note that $f_*R_1=f_*R_2$ implies $R_1=R_2$ 
for $R_1,R_2\in \mathcal{S}$). 
Moreover since we may assume that $f_*z$ is a rational 1-cycle by the induction hypothesis (the rationality of the geometrically
extremal rays of $\NEbar (Z)$), combining the short exact sequence 
above with the fact $K_Y\cdot z=0$ and $K_Y \cdot F\in\Q_{<0}$, we can conclude that we may take $z$ as a rational 1-cycle. 
Use Theorem \ref{ans. of c_2 conjecture}(i) again, we have that the set $\{f_R\}_{R\in \mathcal{S}}$ is finite and in particular
$\mathcal{S}$ is finite. This completes the proof.

(ii)
We may assume that $\varphi$ is primitive. Put $\textup{B}_{\Delta}:=\bigl\{ \alpha \in \Aut X \bigm| \alpha _* \Delta \subset \varphi ^*\mathcal{A}^e (Y) \bigr \}$ for a codimension 1 face $\Delta$ of $\Pi$ and 
$\textup{B}:=\coprod_{\Delta\subset\Pi}\textup{B}_{\Delta}$, where $\Delta$ runs through every codimension 1 face of $\Pi$.
Then we have 
\[\varphi ^*\sAbar (Y)=\overline{\varphi ^* \mathcal{A}^e (Y)}=\overline{\bigcup_{\alpha \in \textup{B}} (\alpha _* \Pi \cap \varphi ^* \mathcal{A}^e (Y))}.\] 
Here we take the closure in the relative topology of the real vector subspace $\Span{\varphi ^* \mathcal{A}^e (Y)}\subset N^1(X)$.
Hence it is enough to prove that $\textup{B}_{\Delta}$ is a finite set for every $\Delta$. Fix a codimension 1 face $\Delta$ 
such that $\textup{B}_{\Delta}\ne\emptyset$. Replace $\Pi$ with $\alpha _*\Pi$ for some $\alpha\in \Aut (X)$ if necessary,
then we may assume that $\Delta \subset \varphi ^*\mathcal{A}^e (Y)$. First we look for classes of ample divisors on $Y$ on 
which $\varphi _* c_2$ takes minimum value and whose pull back on $X$ belongs to $\Delta$. Since 
$\varphi ^*\sAbar (Y)\subset {c_2}_{>0}$, there are only finitely many such and by adding
these together and pulling it back on $X$, we get a nef and big divisor $H$ on $X$. Of course $[H]\in \Delta$ by the definition.
Note that the set $\{[\alpha _*H]\}_{\alpha\in\textup{B}_{\Delta}}$ is finite and so put this by 
$\{[{\alpha _1}_* H],\dots ,[{\alpha _n}_* H]\}$,
where $\alpha_i\in\textup{B}_{\Delta}$. It is straightforward to see that $\textup{B}_{\Delta}=\coprod_{i=1}^n \alpha _i\cdot \Aut(X,H)$. Therefore
we know that $\textup{B}_{\Delta}$ is a finite set by Lemma \ref{nef&big} and the proof is done.    
\end{proof}


\section{The structure of certain C--Y 3-folds with infinitely many divisorial contractions}\label{section structure}
The main results of this section are Theorem \ref{structure} and Corollary \ref{I_{c_2=0}}. We use the following notation and terminology.
\renewcommand{\labelenumi}{(\roman{enumi})}
\begin{enumerate}

\item Let $X$ be a normal projective variety such that $\mathcal{O}_X(K_X)\simeq \mathcal{O}_X$. We denote by $\omega _X$ a generator of $H^0(X,\mathcal{O}_X(K_X))$. A finite automorphism group $G$ is called $\emph{Gorenstein}$ if $g^* \omega _X=\omega _X$ for all $g\in G$. 

\item Suppose we have a faithful finite group action $G$ on a variety $X$. Put 

$X^g:=\bigl\{ x\in X \bigm| g(x)=x \bigr\}$ for $g\in G$; $\quad X^{[G]}:=\bigcup _{g\in G\backslash \{1\}}X^g$.

\item Put $\zeta _n:=\exp (2\pi i /n)$, the primitive $n$-th root of unity in $\C$. Denote by $E_{\zeta}$ the elliptic curve whose period is 
$\zeta$ in the upper half plane. Let us recall the following pairs of an abelian 3-fold and its specific Gorenstein automorphism group: the 
pair $(A_3,g_3)$, where $A_3$ is the triple product of $E_{\zeta _3}$ and $g_3$ is its automorphism $\diag (\zeta _3,\zeta _3,\zeta _3)$ and the 
pair $(A_7,g_7)$ is the Jacobian 3-fold of the Klein quartic curve $C=(x_0x_1^3+x_1x_2^3+x_2x_0^3=0)\subset\PP ^2$ and $g_7$ is the automorphism 
of $A_7$ induced by the automorphism of $C$ given by $[x_0:x_1:x_2]\mapsto [\zeta _7x_0:\zeta _7 ^2x_1:\zeta _7 ^4x_2]$. We call $(A_3,g_3)$ a 
\emph{Calabi pair} and $(A_7,g_7)$ a \emph{Klein pair}.    
\end{enumerate}


\begin{defn} 
Let $W$ be a normal projective surface over $\C$ with at most klt singularities. We call $W$ a \emph{log Enriques surface} if $h^1(\mathcal{O}_W)=0$, $m K_W=0$ for some positive integer $m$. We call the integer $I(W):=\min \bigl \{ m\in \Z_{>0} \bigm | mK_W=0 \bigr \}$ the global canonical index of $W$. 
\end{defn}

\noindent
We construct C--Y 3-folds with infinitely many birational contractions from certain log Enriques surfaces in Section \ref{construction}.


\begin{defn} Let $\varphi :X\to Y$ be a contraction from a C--Y 3-fold $X$ and a divisor $L$ on $X$ the pull back of an ample divisor on $Y$. 
We call $\varphi$ a \emph{$c_2$-contraction} if $L\cdot c_2 =0$. For example, a fibration $\varphi :X\to \PP ^1$ is a $c_2$-contraction if and 
only if the general fiber is an abelian surface. Moreover for an elliptic fibration $\varphi :X\to W$, it is a $c_2$-contraction if and only if 
$W$ is a log Enriques surface by \cite{Oguiso93} (we do not have to assume there that $X$ is simply connected). There exists a unique 
$c_2$-contraction $\varphi _0 :X\to Y_0$ such that every $c_2$-contraction $\varphi :X\to Y$ on $X$ factors through $\varphi _0$ (see 
Lemma-Definition (4.1), \cite{OS}). We call $\varphi _0$ the \emph{maximal $c_2$-contraction}.
\end{defn}  

\noindent
We have the beautiful classification of C--Y 3-folds which admit either a birational $c_2$-contraction or an elliptic $c_2$-contraction, due to K. Oguiso (see \cite{OS}). It plays an important role to prove Theorem \ref{structure}. The following result is coarser than the Oguiso's original classification.

  %
  %
  %
  
\begin{thm}[Oguiso]\label{O-S}
  
\renewcommand{\labelenumi}{(\roman{enumi})}
\begin{enumerate} 
\item Let $\varphi :X\to Y$ be a non-isomorphic birational $c_2$-contraction. Then $\varphi$ is isomorphic to either one of the following:
\renewcommand{\labelenumi}{(\roman{enumi})}
\begin{enumerate} 
\item The unique crepant resolution $\Phi_{7} : X_{7} \rightarrow \bar{X}_7:=A_7/\Span{g_7}$ of $\bar{X}_7$, where $(A_7,g_7)$ is the Klein pair. 
\item The unique crepant resolution $\Phi_{3} : X_{3} \rightarrow \bar{X}_3:=A_3/\Span{g_3}$ of $\bar{X}_3$, where $(A_3,g_3)$ is the Calabi pair.
\item The unique crepant resolution $\Phi_{3,i} : X_{3,i} \rightarrow \bar{X}_{3,i}$ of $\bar{X}_{3,i}$, $(i=1,2)$, where $\bar{X}_{3,i}$ is an \'etale quotient of $\bar{X}_3$.
\end{enumerate}

\item Let $\varphi :X\to W$ be an elliptic $c_2$-contraction. Then $\varphi$ is isomorphic to either one of the following:

\renewcommand{\labelenumi}{(\roman{enumi})}
\begin{enumerate} 
\item
One of the relatively minimal models over 
$W_3$ of 
$$p_{12} : X_3 \overset{\Phi_{3}}\longrightarrow 
\overline{X}_3 
\overset{\overline{p}}\longrightarrow W_3,$$ 
where $\Phi_{3} : X_{3} \rightarrow \bar{X}_3$ is as above and $\overline{p}$ is an elliptic fibration on $\overline{X}_3$.

\item 
An elliptic fiber space structure on an \'etale quotient of an abelian 3-fold. 

\item  
One of the relatively minimal models over $W_{3,1}$ of 
$$\kappa_{3,1} : X_{3,1} \overset{\Phi_{3,1}}\longrightarrow 
\overline{X}_{3,1}  
\overset{\overline{\kappa}}\longrightarrow W_{3,1},$$ 
where $\Phi_{3,1} : X_{3,1} \rightarrow \bar{X}_{3,1}$ is as above and $\overline{\kappa}$ is an elliptic fibration on $\overline{X}_{3,1}$.

\item
One of the relatively minimal models over $S/G$ of 
$$\psi :Y \overset{\nu}\rightarrow 
(S\times E)/G \overset{\mu}\rightarrow S/G,$$ 
where $S$ is a normal K3 surface (namely its minimal resolution is a smooth K3 surface), $E$ is an elliptic curve, 
$G$ is a finite Gorenstein automorphism group of $S \times E$ 
whose element is of the form $(g_{S}, g_{E}) \in \Aut S 
\times \Aut E$ and $\nu$ is a crepant resolution 
of $(S\times E)/G$. Slightly more precisely, $G$ is of the form 
$G = H \rtimes \Span{a}$, where $H$ is a commutative group 
consisting of elements like $h = (h_{S}, h_{E})$
such that $\text{ord}(h_{S}) = \text{ord}(h_{E}) = \text{ord}(h)$ and 
$h_E$ is a translation, furthermore the generator 
$a$ of $\Span{a}$ is the element of the form $(a_{S}, \zeta_{I(W)}^{-1})$ such that  
$a_{S}^{*}\omega_{S} = \zeta_{I(W)}\omega_{S}$. Moreover $I(W) \in \{2, 3, 4, 6 \}$.   
\end{enumerate}
\end{enumerate}
\end{thm}

For a contraction $\varphi :X\to Y$ on a C--Y 3-fold $X$, we define $M(\varphi):=\bigl\{i \in I\bigm| E_i\cdot C=0$ for all curves $C$ such that $\varphi (C)$ is a point $\bigr \}$.

%
%
%

\begin{lem}\label{special3}
\renewcommand{\labelenumi}{(\roman{enumi})}
\begin{enumerate} 
\item 
Let $\varphi :X\to Y$ be a primitive birational contraction on a C--Y 3-fold $X$. Denote the extremal ray corresponding to $\varphi$ by $R$. Then the set 
\[
  \mbox{$L(\varphi):=\bigl\{i \in I\bigm| R\subset V_i$ and $\varphi (E_i)$ is a $\Q$-Cartier divisor on $Y \bigr \}$}
  \]
is finite. 

\item 
Let $\varphi :X\to Y$ be a (not necessarily primitive) birational contraction on a C--Y 3-fold $X$. The set
\begin{align*}
  \overline{M(\varphi)} & :=\bigl\{ i\in M(\varphi) \bigm| E_i \cap \Exc (\varphi) \ne \emptyset \bigr \} \\ & =\bigl\{ i\in I \bigm| E_i \cap \Exc(\varphi) \ne \emptyset 
  \text{ and } E_i= 0 \text{ in } N^1(X/Y) \bigr\}
\end{align*}
is finite. 

\item 
Suppose that we have the following diagram:
\[ \xymatrix{ X \ar[dr]_{\varphi} \ar@{-->}[rr] ^{\Phi\quad} & & Y \ar[dl]^{\psi} \\
   & W & ,}\]
where $\varphi$, $\psi$ are contractions on C--Y 3-folds $X$, $Y$ and $\Phi$ is a birational map over $W$. Then for general $i\in M(\varphi)$, $E_i$ is contained in the isomorphic locus of $\Phi$. In particular, $|M(\varphi)|=\infty$ is equivalent to $|M(\psi)|=\infty$.
\end{enumerate}
\end{lem}  

\begin{proof}
(i) Assume that $L(\varphi)$ is infinite. We can take $1,2\in L(\varphi)$ such that $E_1\cap E_2 \ne \emptyset$. Since $R\subset V_1\cap V_2$, the class of 1-cycle $[E_1\cdot E_2]$ belongs to $R$ and so $\dim \varphi (E_1\cap E_2)=0$. Hence $\dim \varphi (E_1)\cap \varphi (E_2)=0$. This is a contradiction because $\varphi (E_1)$ and $\varphi (E_2)$ are $\Q$-Cartier divisors.

(ii) Let $R_1,\dots ,R_n$ be the generators of the cone $\overline{NE}(X/Y)$, namely extremal rays, and consider that $\psi_k$ is the extremal contraction corresponding to $R_k$. It is enough to check that $\overline{M(\varphi)} \subset \bigcup _{k=1}^{n}L(\psi_k)$. Pick $0 \in \overline{M(\varphi)}$. Then there exist an integer $k$ and an irreducible curve $C$ such that $C\subset E_0$ and $[C]\in R_k$. Thus $R_k\subset V_0$. Now since $\psi_k (E_i)$ is a Cartier divisor for $i\in M(\varphi)$, we obtain the statement.
 
(iii) Note that $\Phi$ is a composition of flops over $W$. Apply (ii) for each flopping contraction, then we obtain the statement. 
\end{proof}

%
%
%

\begin{lem}\label{typeA}
We use the notation in Theorem \ref{O-S}. Neither $X_7$, $X_3$, $X_{3,1}$ nor $X_{3,2}$ admits infinitely many contractions of type III.
\end{lem}

\begin{proof}
Let $\Phi _3$ be the unique crepant resolution of $\bar{X}_3$. $\Phi _3$ is a composition of birational contractions of type II (cf. \cite{Oguiso96}). Pick $i\in I_{X_3}$, if any. Then $\Phi _3(E_i)\cap \Sing \bar{X}_3 \ne \emptyset$ because $\bar{X}_3$ is a quotient of an abelian 3-fold. Since $\Sing \bar{X}_3=\Phi_3(\Exc (\Phi_3))$, we have $E_i\cap \Exc (\Phi)\ne\emptyset$, which implies $i\in L(\psi )$ for some contraction $\psi$ of type II.
Hence if $I_{X_3}$ is infinite, there exists a birational contraction $\psi$ of type II on $X_3$ such that $L(\psi)$ is infinite. This is absurd.
In the cases of $X_{3,1}$ and $X_{3,2}$, the same proof as above works, since $\bar{X}_{3,1}$, $\bar{X}_{3,2}$ are \'{e}tale quotients of
$\bar{X_3}$. Next let $\Phi _7$ be the unique crepant resolution of $\bar{X_7}$. Then $\Exc (\Phi _7)=E_1\cup E_2\cup E_3$, each $E_j$ is
a Hirzeburch surface of degree 2 and these divisors are crossing normally each other along the negative sections (cf. \cite{Oguiso96}) 
(thus $v_a\in\R _{\ge 0}[F_b]$, $v_b\in\R _{\ge 0}[F_c]$, $v_c\in\R _{\ge 0}[F_a]$ for some $a,b,c$ with $\{a,b,c\}=\{1,2,3\}$).
Because $\bar{X}_7$ is a quotient of an abelian 3-fold, $E_i\cap (E_1\cup E_2 \cup E_3)\ne\emptyset$ for all $i\in I_{X_7}$.
Furthermore if $E_i$ intersects $E_a$ and if $v_i\notin \R_{\ge 0}[F_b]$, $v_i\in \R_{\ge 0}[F_a]$, since 
$v_i \in V_a\cap E_b ^{\perp}$. So in this case $E_i$ intersects $E_a$ and $E_c$, does not intersect $E_b$.
By this way, we know that every $E_i$ intersects precisely two of $E_1, E_2$ and $E_3$. Assuming that $I_{X_7}$ is infinite, we can find
a divisorial contraction $\psi$ which contracts either $E_1, E_2$ or $E_3$, such that $L(\psi)$ is infinite.
So we obtain a contradiction. 
\end{proof} 

  %
  %
  %
  
\begin{thm}\label{structure}
Assume that $I_{c_2=0}(=I_{X,c_2=0}):=\bigl\{i \in I_X\bigm| E_i\cdot c_2=0 \bigr \}$ is infinite. Then the following hold.
\renewcommand{\labelenumi}{(\roman{enumi})}
\begin{enumerate} 

\item 
We have an elliptic $c_2$-contraction $\varphi :X\to W$ and $\varphi$ fits in the case of (ii)$(d)$ in Theorem \ref{O-S}, that is, we have the following diagram: 
\[ \xymatrix{ X \ar[dr]_{\varphi} \ar@{-->}[rr] ^{\Phi\quad}& & Y \ar[dl]^{\psi} 
\ar[d]^{\nu}\\
   & W\cong S/G & (S\times E)/G \ar[l]^{\mu},}\]
where $Y$, $S$, $E$, $G$, $\psi$, $\nu$ and $\mu$ are given there. Let $r: S\times E\to (S\times E)/G$ be the quotient morphism. Then the normal K3 surface $S$ contains infinitely many smooth rational curves $\{l\}$ such that
\renewcommand{\labelenumi}{(\roman{enumi})}
\begin{enumerate} 
\item $r(l\times E) \cap \Sing (S\times E)/G =\emptyset$, and
\item $\bigcup _{g\in G} g\cdot l$ is contractible at the same time by a birational contraction on $S$. 
\end{enumerate}

\item Let $\Phi$ denote the birational map between $X$ and $Y$ over $W$ in (i). Then for general $i\in I_{c_2=0}$, $E_i$ is contained in the isomorphic locus of the birational map $\nu \circ \Phi$ and $E_i =r(l\times E)$ under this isomorphism for some smooth rational curve $l$ on $S$ satisfying (a) and (b) in (i).
\end{enumerate}
\end{thm}

\begin{proof}
(i) Let us denote by $\varphi :X\to W$ the maximal $c_2$-contraction (a priori $W$ may be a point).   

\begin{cla}\label{claim I_{c_2=0}}
For a general $i\in I_{c_2=0}$, $i\in M(\varphi)$.
\end{cla}

\begin{proof}
If not, by Proposition \ref{interesting} we can take $1,2,3 \in I_{c_2=0} \backslash M(\varphi)$ such that some multiple of $E_1+E_2+E_3$ determines a $c_2$-contraction, which factors through $\varphi$. By the choice of $1,2,3$, there exists one of the elements $1,2,3$, say $1$, and there exists an irreducible curve $C$ on $X$ such that $\varphi (C)$ is a point and $E_1\cdot C>0$. By the proof of \ref{interesting} we can pick $4,5\in I_{c_2=0} \backslash M(\varphi)$, different from $1,2,3$, such that some multiple of $E_1+E_4+E_5$ determines a $c_2$-contraction, which factors through $\varphi$. Thus there exists one of the elements $4,5$, say $4$, such that $E_4\cdot C<0$. By the same procedure, we have infinitely many elements $i \in I_{c_2=0} \backslash M(\varphi)$ such that $E_i\cdot C<0$. This is a contradiction with \ref{interesting}.  
\end{proof}

\noindent When $\dim W=1$, at most finitely many $E_i$ $(i\in I)$ are contracted to a point on $W$ by $\varphi$, so $M(\varphi)$ is finite.
Hence we have $\dim W\ge 2$. If $\varphi$ is isomorphic, $\sAbar (X)\subset c_2 ^{\perp}$ and in particular $c_2 =0$. In this case, 
$X$ is an \'etale quotient of an Abelian 3-fold by \cite{Kobayashi} and it never admits birational contractions.
Combining Theorem \ref{O-S} with Lemma \ref{special3}(iii) and Lemma \ref{typeA}, we know that $\varphi$ fits in the case (ii)$(d)$ of \ref{O-S} and $|M(\psi)|=\infty$. Furthermore $|M(\psi)|=\infty$ implies
that the set $\bigl\{i \in I_{(S\times E)/G} \bigm| E_i\cap \Sing (S\times E)/G =\emptyset \bigr \}$ is infinite by \ref{special3}(ii).
Here we use the equality $\Sing (S\times E)/G=\nu (\Exc (\nu))$. Note that every primitive birational contraction on $S\times E$ is the form
as $f\times id_E$, where $f$ is a contraction of a single smooth rational curve on $S$. Thus we have the conditions $(a)$ and $(b)$.

(ii) This follows from \ref{special3}(ii) and \ref{special3}(iii).  
\end{proof}

  
\begin{rem}\label{converse}
\renewcommand{\labelenumi}{(\roman{enumi})}
\begin{enumerate}
\item Assume that Theorem \ref{structure}(i) holds. Then we have an infinite set $\bigl\{i \in M(\mu) \bigm| E_i\cap \Sing (S\times E)/G =\emptyset \bigr \}$. Using Lemma \ref{special3}(iii), we know that $I_{X,c_2=0}$ is infinite. Namely \ref{structure}(i) is a characterization of C--Y 3-folds $X$ with $|I_{X,c_2=0}|=\infty$. 
\item  Because $(\Sing S\times E)\cup (S\times E)^{[G]}=r^{-1}\Sing (S\times E)/G$ by the purity of branch locus, the condition $(a)$ in \ref{structure}(i) is equivalent to the condition
\[(a)' \quad(l\times E)\cap ((\Sing S\times E)\cup (S\times E)^{[G]})=\emptyset.\]
\end{enumerate}
\end{rem}

%
%
  
\begin{cor}\label{I_{c_2=0}}
The set $I_{c_2=0}$ is finite up to $\Aut X$.
\end{cor}

\begin{proof} We may assume that $I_{c_2=0}$ is infinite. Now $X$ is birational to $(S\times E)/G$ via $\nu \circ \Phi$ as in Theorem \ref{structure}. Consider the minimal resolution $S'\to S$. We may assume that $Y$ is obtained as a crepant resolution $\nu ':Y\to (S'\times E)/G$, that is, $\nu$ factors through $\nu '$. The existence of $\nu '$ is guaranteed by \cite{Roan}. By \ref{structure}(ii) and Claim \ref{claim I_{c_2=0}}, for general $i \in I_{c_2=0}$, $E_i$ is contained in the isomorphic locus of $\nu '\circ \Phi$ and $E_i$ is isomorphic to the image on $(S'\times E)/G$ of $l\times E$ for some smooth rational curve $l$ on $S'$. On the other hand, the set $I_{(S'\times E)/G}$ is finite up to $\Aut(S'\times E)/G$ by Theorem (2.23) in \cite{OS} (note that the proof of Theorem (2.23) in \cite{OS} works even if $G$ does not act on $S'\times E$ freely). Therefore the set $I_{c_2=0}$ is finite up to $\Bir X$. By the proof of Lemma (1.15) in \cite{K97}, the set $I_{c_2=0}$ is
  finite up to $\Aut X$. 
\end{proof}

\noindent 
As we mention in the Introduction, the following problem seems worthwhile to think about.

%
%

\begin{prob}\label{existence of c_2}
Assume that $\Aut X$ is infinite and its Picard number $\rho (X)$ is sufficiently large. Then does $X$ admit a nontrivial $c_2$-contraction?
\end{prob}

\noindent
Conjecture \ref{c_2 conjecture} says that if $\Aut X$ is infinite the shape of $\sAbar (X)$ is 
complicated near $\sAbar (X)\cap c_2 ^{\perp}$. We expect that this ``complexity'' produces a rational point on
$\sAbar (X)\cap c_2 ^{\perp}\backslash \{0\}$ and some multiple of the divisor corresponding to the rational point defines a 
$c_2$-contraction.  
In fact when we study the structure of C--Y 3-folds $X$ with $|I_{c_2 =0}|=\infty$ in Theorem \ref{structure}, we showed 
the existence of an elliptic $c_2$-contraction on $X$ by Proposition \ref{interesting}. 



\section{Construction of C--Y 3-folds with infinitely many birational contractions}\label{construction}
The aim of this section is to give construction of C--Y 3-folds with infinitely many birational contractions of type I or III from certain log Enriques surfaces. First of all, given a log Enriques surface $W$ with $I(W)\in\{2,3,4,6\}$, we construct a C--Y 3-fold $X$ with a $c_2$-contraction $\varphi: X\to W$. Let $q:S\to W$ be the global canonical cover and denote by $G=\Span{a}(\cong \Z /I(W)\Z)$ the Galois group of $q$. The $S$ may be an abelian surface in general but here we assume that $S$ is a normal K3 surface (this assumption is satisfied, for example, if $W$ contains a \emph{contractible} smooth rational curve. Here a curve $m$ on $W$ is said \emph{contractible} if it is contracted by a birational contraction and this is equivalent to $m^2<0$). Let $E$ be an elliptic curve such that $E$ has an automorphism of order $I(W)$ which fixes the origin. Suppose that the generator $a$ of $G$ satisfies that $a^{*}\omega_{S} = \zeta_{I(W)}\omega_{S}$. Then define the action of $
 a$ on $E$ as $a(x)=\zeta_{I(W)}^{-1}x$ for $x\in E$. Then $G$ gives a Gorenstein action on $S\times E$. Take the minimal resolution $S'\to S$, then $G$ acts on $S'$ and we know that $(S'\times E)/G$ is a C--Y model. By \cite{Roan} there exists a crepant resolution $\nu ':X\to (S'\times E)/G$. Of course this $X$ is a C--Y 3-fold and $\varphi :X \rightarrow 
(S'\times E)/G \rightarrow (S\times E)/G  \rightarrow S/G=W$ is an elliptic $c_2$-contraction.

 For a log Enriques surface $W$, let us denote by $\Sigma _W$ the locus of klt points on $W$ which are neither RDP's nor smooth points.

%

\begin{prop}\label{type I and III}
Let $\varphi :X \overset{\nu '}\rightarrow 
(S'\times E)/G  \overset{\mu}\rightarrow S/G=W$ be as is constructed from $W$ above. Suppose that there exists a contractible smooth rational curve $m$ on $W$.
\renewcommand{\labelenumi}{(\roman{enumi})}
\begin{enumerate}
\item Assume that $m\cap\Sigma _W=\emptyset$. Then there exists a contraction of type III on $X$ contracting a prime divisor $D_0$ such that $\varphi (D_0)=m$.

\item Assume that $m\cap\Sigma _W\ne\emptyset$. Then there exists a contraction of type I on $X$ contracting an irreducible curve $m_0$ such that $\varphi (m_0)=m$.
\end{enumerate}
\end{prop}

\begin{proof}
Let $r': S'\times E\to (S'\times E)/G$ be the quotient morphism. Moreover let $l$ be an irreducible component of $q^{-1}m$
and denote by $l'$ the strict transform of $l$ on $S'$. Put $D:=r'(l'\times E)$. In the first case, because $l'\cap {S'}^{[G]}=\emptyset$
we know that $D\cap \Sing (S'\times E)/G =\emptyset$. Furthermore since $m$ is contractible on $W$, $\bigcup _{g\in G} g\cdot l'$ is
contractible on $S'$ and in particular, $D$ is contractible by a birational contraction of type III on $(S'\times E)/G$. Hence ${\nu '}_*^{-1}D$
gives a desired divisor $D_0$. In the second case, we have $(l\times E)\cap (S\times E)^{[G]}\ne\emptyset$ (we prove the contraposition of this
in the proof of Proposition \ref{log Enriques} below) and $D$ is an exceptional divisor of a contraction of type III, since
$\bigcup _{g\in G} g\cdot l'$ is contractible on $S'$. Moreover $D$ contains a point $y\in r'((S'\times E)^{[G]})$ such that $y$ is over
a point in $m\cap \Sigma_W$ by the morphism $\mu$. Note that $\dim (S'\times E)^{[G]}\cap (l'\times E)=0$. Because the problem is local, we
may assume that $\{y\}=(\Sing (S'\times E)/G)\cap D$. Let 
\[ X=:X_0\overset{\psi_1}\longrightarrow X_1 \cdots \overset{\psi_n}\longrightarrow X_n:=(S'\times E)/G\]
 be a primitive decomposition of $\nu '$ and let us denote by $m_n$ the unique irreducible curve passing through $y$, of the form
$r'(l'\times \{z\})$, where $z$ is a point in $E^{[G]}$. Suppose that $D_i$ (resp. $m_i$) stands for the strict transform of $D$
(resp. $m_n$) on $X_i$. Let $V$ be an irreducible component of ${\nu '}^{-1}y$ such that $V\cap D_0\ne \emptyset$. When $\dim V=2$, we have
$\dim V\cap D_0=1$. If every component $V$ such that $V\cap D_0\ne \emptyset$ is $1$-dimensional, the equality ${\nu '}^* D\cdot V=0$ implies
that $V\subset D_0$, hence $\dim V\cap D_0 =1$ (note that $D_0$ is not contractible any more by a divisorial contraction on $X$, since the
dimension of the image of the map $N_1(D_0)\to N_1(X)$ is more than 2 (cf. Fact (iii))). Therefore there exists an integer $k\ge 1$ such that
$\dim \psi^{-1}_{k+1}\cdots \psi^{-1}_{n}y\cap D_k=0$ and $\dim \psi^{-1}_{k}\cdots \psi^{-1}_{n}y\cap D_{k-1}=1$. The following claim comes 
from the general theory and we leave the proof to the reader, since it is an easy exercise.   

\begin{cla}
Let $f:X\to Y$, $g:Y\to Z$ be primitive birational contractions between C--Y models. Suppose that the strict transforms $f^{-1}_*l$ of all curves $l$ contracted by $g$ are numerically proportional. Then if $g$ is of type I (resp. of type III), there exists a contraction $f'$ of type I (resp. of type III) over $Z$ such that $f^{-1}_*l$ are contracted by $f'$.   
\end{cla} 

\noindent
We apply the claim repeatedly and then we have a contraction of type III on $X_k$, $\psi: X_k\to Z$, such that $\Exc (\psi)=D_k$. Let $C_{k-1}$ be an irreducible curve on $X_{k-1}$ such that $C_{k-1}\subset \psi^{-1}_{k}\cdots \psi^{-1}_{n}y\cap D_{k-1}$. Then we know that $\overline{NE}(X_{k-1}/Z)$ is generated by $\R_{\ge 0}[C_{k-1}]$ and $\R_{\ge 0}[m_{k-1}]$. The latter extremal ray determines a contraction of type I on $X_{k-1}$ and using the claim again, we obtain a contraction of type I on $X$ whose exceptional set consists of $m_0$. 
\end{proof}

\noindent
Consider a log Enriques surface $W$ with $I(W)\in\{2,3,4,6\}$ such that $W$ contains infinitely many contractible smooth rational curves. Then by Proposition \ref{type I and III}, we can construct a C--Y 3-fold $X$ with infinitely many birational contractions of type I or type III. 

%
%

\begin{ex}
\renewcommand{\labelenumi}{(\roman{enumi})}
\begin{enumerate}
\item
See the nice survey, \cite{Kon}, by S. Kond\={o} and its references for the details of the following. Due to E. Horikawa we know that
the moduli space $\mathcal{M}$ of Enriques surfaces is $10$-dimensional. The moduli space $\mathcal{N}$ of Enriques surfaces which contains
at least one smooth rational curve is an irreducible subvariety of codimension $1$ in $\mathcal{M}$. Enriques surfaces whose automorphism
group is finite are classified by S. Kond\={o} and the moduli of them consists of seven families $\{ \mathcal{F}_i \} _{i=1} ^7$ and each
family is at most $1$-dimensional. On the other hand for Enriques surfaces $W$, $\Aut W$ is finite if and only if $W$ contains at least one
but at most finitely many smooth rational curves. Consequently there exists the $9$-dimensional moduli space,
$\mathcal{N}\backslash \bigcup_{i=1}^7{\mathcal{F}_i}$, whose elements are Enriques surfaces which contain infinitely many smooth rational curves. 
\item
Let $E_1$, $E_2$ be elliptic curves which are not mutually isogenous and $S'$ the Kummer surface associated to the abelian surface $E_1\times E_2$.
Consider the involution $a$ on $S'$ induced by the involution $(x,y)\mapsto (x,-y)$ on $E_1\times E_2$. Let $\{F_i\}_{i=1}^4$
(resp. $\{F' _i\}_{i=1}^4$) be the smooth rational curves on $E_1\times E_2/(-1)$ which are
the images of $\{x\}\times E_2$ (resp. $E_1 \times \{y\}$) by the natural map $E_1\times E_2\to E_1\times E_2/(-1)$,
where $x\in E_1$ (resp. $y\in E_2$) is a point of order 2. Then the fixed locus ${S'}^{a}$ consists of the eight, disjoint smooth rational
curves $f_* ^{-1}F_i$, $f_* ^{-1}F' _i$, where $f$ is the minimal resolution of $E_1\times E_2/(-1)$. Because the every generator of the
Picard group of $S'$ is fixed by the involution $a$, every smooth rational curve $l'$ is also fixed, that is, $a\cdot l'=l'$. Contract the
eight smooth rational curves $f_* ^{-1}F_i$, $f_* ^{-1}F' _i$ on $S'$ and we get a normal K3 surface $S$ with eight $\rm{A_1}$-singularities.
The group action of $\Span{a}$ on $S'$ descends to the group action on $S$ and let us use the same letter $\Span{a}$ for this action. Then we
obtain a log Enriques surface $W:=S/\Span{a}$ which contains infinitely many contractible smooth rational curves $\{m\}$ such that
$m\cap \Sigma _W\ne \emptyset$. Here we use the fact that every Kummer surface has the infinite automorphism group and so in particular,
it contains infinitely many smooth rational curves. 
\end{enumerate}
\end{ex}

\noindent 
I do not know any example of \emph{rational} log Enriques surface $W$ which contains infinitely many smooth rational curves 
$\{m\}$ such that $m\cap \Sigma _W=\emptyset$.\footnote{If a log Enriques surface $W$ satisfies such conditions, the minimal 
resolution of $W$ contains infinitely many $-2$ curves. I found an example of a smooth rational surface containing 
infinitely many $-2$ curves but unfortunately my surface is not the minimal resolution of log 
Enriques surface.}

 The following statement is the converse of Proposition \ref{type I and III} .
%
 
\begin{prop}\label{log Enriques}
Suppose the conditions in Theorem \ref{structure}(i) hold. Then the log Enriques surface $W\cong S/G$ contains infinitely many contractible
smooth rational curves $\{m\}$ such that $m=\varphi (E_i)$ and $m\cap \Sigma_W =\emptyset$.
\end{prop}

\begin{proof} 
Because $G = H \rtimes \Span{a}$ as is in (ii)$(d)$ in Theorem \ref{O-S}, we can decompose the quotient morphism 
$S\to W$ as follows:
\[ \xymatrix{S \ar[r]^{p\qquad} & T:=S/H \ar[r]^{q\quad} & S/G=T/\Span{a}\cong W}.\]
Note that $T$ is a normal K3 surface, for $H$ is a Gorenstein group acting on $S$ (and notice that $H$ was trivial in the argument before
Proposition \ref{type I and III}). In particular, $T$ has at most RDP's. 

\begin{cla}
$l\cap S^{h\cdot a^i}=\emptyset$ for all $h\in H$, all $i\ne 0$ modulo $I(W)$.
\end{cla}

\begin{proof} 
The condition Remark \ref{converse}$(a)'$ implies that $(l\times E)\cap (S\times E)^{[G]}=\emptyset$. Therefore if $E^{h\cdot a^i}\ne \emptyset$
for all $h\in H$, all $i\ne 0$ modulo $I(W)$, we know that $l\cap S^{h\cdot a^i}=\emptyset$. In fact this hypothesis is true, since the morphism
$id_E-a^i$ on $E$ is surjective.
\end{proof}

\noindent
It is straightforward to see that 
\[p^{-1}T^{a^i}=\bigcup _{h\in H}S^{h\cdot a^i} \text{ for all }i.\] 
Thus we have $p(l)\cap T^{[\Span{a}]}= \emptyset$. On the other hand because $W\backslash q(T^{[\Span{a}]})$ has at most RDP's,
$q\circ p(l)\cap \Sigma_W =\emptyset$. Since $q\circ p(l)$ is contractible by an extremal contraction on $W$, $q\circ p(l)\cong \PP^1$.
\end{proof}

\noindent
In summary, for a given C--Y 3-fold $X$ with $|I_{c_2=0}|=\infty$ there exists an elliptic $c_2$-contraction $\varphi :X\to W$. Here $W$ is
a log Enriques surface with $I(W)\in\{2,3,4,6\}$ which contains infinitely many smooth rational curves $\{m\}$ such that
$m\cap \Sigma_W =\emptyset$ and $m=\varphi (E_i)$ for some $i\in I_{c_2=0}$. Conversely, for a given log Enriques surface $W$ with
$I(W)\in\{2,3,4,6\}$ which contains infinitely many smooth rational curves $\{m\}$ such that $m\cap \Sigma_W =\emptyset$, there exists
a C--Y 3-fold $X$ with $|I_{c_2=0}|=\infty$ which admits an elliptic $c_2$-contraction $\varphi :X\to W$.

\indent

\textsc{Research Institute for Mathematical Sciences, Kyoto University, Kyoto, 606-8502, Japan} \\
\texttt{hokuto@kurims.kyoto-u.ac.jp}

\end{document}